\documentclass[12pt]{article}
\usepackage{amsmath,amssymb,amsthm}
\usepackage{cases}
\usepackage{amsfonts}
\usepackage{color,xcolor,cite}
\usepackage[left=2.0cm,right=2.0cm,top=2.0cm,bottom=2.0cm]{geometry}
\usepackage[colorlinks,citecolor=blue,urlcolor=blue]{hyperref}

\newtheorem{theorem}{Theorem}[section]

\newtheorem{lemma}{Lemma}[section]
\newtheorem{proposition}{Proposition}[section]

\newtheorem{definition}{Definition}[section]

\newcommand{\bal}{\begin{align}}
\newcommand{\bbal}{\begin{align*}}
\newcommand{\beq}{\begin{equation}}
\newcommand{\eeq}{\end{equation}}
\newcommand{\bca}{\begin{cases}}
\newcommand{\eca}{\end{cases}}

\newcommand{\pa}{\partial}
\newcommand{\fr}{\frac}
\newcommand{\na}{\nabla}
\newcommand{\De}{\Delta}

\newcommand{\cd}{\cdot}

\newcommand{\dd}{\mathrm{d}}

\newcommand{\R}{\mathbb{R}}

\newcommand{\Z}{\mathbb{Z}}

\newcommand{\les}{\lesssim}

\newcommand{\f}{\left}
\newcommand{\g}{\right}

\begin{document}
\title{The failure of H\"{o}lder regularity of solutions for the Euler-Poincar\'{e} equations in Besov spaces }

\author{Guorong Qu$^{1}$ and Min Li$^{2,}$\footnote{E-mail: guorongqu@163.com; limin@jxufe.edu.cn(Corresponding author)}\\
\small $^1$ School of Tourism Data , Guilin Tourism University, Guilin 541006, China\\
\small $^2$ Department of Mathematics, Jiangxi University of Finance and Economics, Nanchang, 330032, China}

\date{\today}

\maketitle\noindent{\hrulefill}

{\bf Abstract:} In this paper, we investigate the continuity of solution to the Euler-Poincar\'{e} equations. We show that the continuity of the solution cannot be improved to the H\"{o}lder continuity. That is, the solution of the Euler-Poincar\'{e} equations with initial data $u_0\in B^s_{p,r}$ belongs to $\mathcal{C}([0,T];B^s_{p,r}(\R^d))$ but not to $\mathcal{C}^\alpha([0,T];B^s_{p,r}(\R^d))$ with any $\alpha\in(0,1)$.


\vskip0mm\noindent{\hrulefill}

\section{Introduction}

In this paper, we are concerned with the Cauchy problem for the
Euler-Poincar\'{e} equations:
\begin{equation}\label{0}
\begin{cases}
\partial_tm+u\cdot \nabla m+\nabla u^T\cd m+(\mathrm{div} u)m=0, \qquad (t,x)\in \R^+\times \R^d,\\
m=(1-\De)u,\qquad (t,x)\in \R^+\times \R^d,\\
u(0,x)=u_0,\qquad x\in \R^d,
\end{cases}
\end{equation}
where $u=(u_1,u_2,\cdots,u_d)$ denotes the velocity of the fluid, $m=(m_1,m_2,\cdots,m_d)$ represents the momentum. To avoid any confusion, in the component-wise, the first equation of  \eqref{0} can be rewrite as
$$\partial_tm_i+\sum^d_{j=1}u_j\partial_{x_j}m_i+\sum^d_{j=1}(\partial_{x_i}u_j)m_j+m_i\sum^d_{j=1}\partial_{x_j}u_j=0,
~~~~~~~~i=1,2,\cdots,d.$$
The Euler-Poincar\'{e} equations \eqref{0} were initially discovered by Holm, Marsden, and Ratiu in \cite{hmr1,hmr2}, serve as a Lagrangian averaged models for modeling and analyzing fluid dynamics, where the Lagrangian is given by the $H^1$ norm of the fluid velocity in $d$-dimensions. These equations have also been analyzed in the context of geodesic motion on the diffeomorphism group as elaborated by Holm and Staley in \cite{hs}. Additionally, in the specific case where $d=2$, the Euler-Poincar\'{e} equations is the same as the averaged template matching equation for computer vision, see references \cite{hma,hrt}.

When $d\geq 2$, the Euler-Poincar\'{e} equations are commonly regarded as multi-dimensional generalization of the Camassa-Holm system. This framework is instrumental in the study and modeling of nonlinear shallow water wave dynamics in Holm et al.'s original paper \cite{hmr1}. Indeed, when $d=1$, the system \eqref{0} simplify to the renowned Camassa-Holm equation
\begin{align}\tag*{(CH)}
m_{t}+um_{x}+2u_{x}m=0, \quad m=u-u_{xx}.
\end{align}
Originally introduced by Camassa and Holm in \cite{ch}, the CH equation is a bi-Hamiltonian model characterizing shallow water wave phenomena. Notably, the CH equation is distinguished by its peakon solutions of the form $Ce^{-|x-Ct|}$, which have generated considerable interest in the physical sciences (see references \cite{c5,t}). The CH equation has been the subject of extensive scholarly inquiry, particularly regarding its well-posedness, weak solution dynamics, and analytic or geometric characteristics. Among the numerous studies, research on local well-posedness and ill-posedness for the CH equation's Cauchy problem is highlighted in \cite{ce2,d2,glmy}. The phenomena of blow-up and the global existence of strong solutions are discussed in \cite{c2,ce2,ce3,ce4}, while the existence of global weak and dissipative solutions has been explored in \cite{bc1,bc2,xz1}, with additional findings cited therein. Furthermore, the non-uniform continuity of the CH equation has been extensively investigated, as evidenced by publications such as \cite{hk,hkm,lyz1,lwyz}.

In terms of the Euler-Poincar\'{e} equations' rigorous analysis, seminal contributions were made by Chae and Liu \cite{cli}, who established the local existence of weak solutions in $W^{2,p}(\R^d),\ p>d$ and the local existence of unique classical solutions in $H^{s}(\R^d),\ s>\frac{d}{2}+3$. Building on this, Yan and Yin \cite{yy} expanded the discourse by discussing the local existence and uniqueness of solutions to \eqref{0} in Besov spaces. Furthermore, Li, Yu, and Zhai \cite{lyzz} provided significant insights by proving that solutions to \eqref{0}, given a broad spectrum of smooth initial data, either exhibit finite-time blow-up or persist globally in time, thereby resolving a question posed by Chae and Liu \cite{cli}. Subsequent research by Luo and Yin led to new findings regarding blow-up in periodic scenarios, leveraging the rotational invariant properties of the equation \cite{luoy}. For a more comprehensive understanding of the Euler-Poincar\'{e} equations, we refer to \cite{luoy,zyl}.

In more recent studies, the focus has shifted towards the continuity properties of the data-to-solution map for Camassa-Holm type equations, a subject of growing interest as indicated in \cite{lyz1}. This aspect is crucial for understanding the well-posedness of these equations. Li, Dai, and Zhu \cite{ldz} demonstrated that the solutions to \eqref{0} exhibit non-uniform continuous dependence on the initial data in $H^{s}(\R^d),s>1+\frac{d}{2}$. This non-uniform continuity result was later extended to the Besov space $B^s_{p,r}(\R^d),s>\max\{1+\frac{d}{2},\frac{3}{2}\}$ in \cite{ldl}. The continuity issue of the data-to-solution map for the Euler-Poincar\'{e} equations in Besov spaces $B^{s}_{p,\infty}(\R^d),\ s>\max\{1+\frac{d}{p},\frac{3}{2}\}$ has been resolved in \cite{lg}, where an ill-posedness result is ultimately established.

As part of the local well-posedness results in Besov spaces of Yan-Yin \cite{yy}, we known that, if $u_0 \in B^s_{p,r}$ with $s > \max\{1+\frac dp,\frac32\},\ (p, r) \in [1, +\infty)^2$ or $s = 1 + \frac{d}{p},\ r=1,\ 1 \leq p < 2d$, then there exists a solution $u \in C([0, T]; B^s_{p,r})$ for the Euler-Poincar\'{e} equations \eqref{0}. In this paper, we are interested in the following

\vspace*{1em}
\noindent\textbf{H\"{o}lder Continuity Problem :} {\it Whether or not the solution $u(t)$ for the Euler-Poincar\'{e} equations \eqref{0} from initial data $u_0 \in B^s_{p,r}$ can belong to $C^\alpha([0, T]; B^s_{p,r})$ with some $\alpha \in (0, 1)$.}
\vspace*{1em}

More specifically, if the initial data \(u_0\) have more regularity such that \(u_0 \in B^{s'}_{p,r}\) for some \(s' > s\), by the interpolation argument, we can deduce that \(u \in C^\alpha([0, T]; B^{s}_{p,r})\) with \(\alpha = s' - s\). In this paper, we will show that there exists initial data \(u_0 \in B^s_{p,r}\) such that the corresponding solution of the Euler-Poincar\'{e} equations \eqref{0} cannot belong to \(C^\alpha([0, T]; B^{s'}_{p,r})\) with any \(\alpha \in (0, 1)\). Namely, we provide a negative response to the question:
$$u_0 \in B^s_{p,r} \stackrel{?}{\implies} u \in C^\alpha([0, T]; B^{s'}_{p,r}) \text{ with } \alpha \in (0, 1).$$
Consequently, this also negates the question posed in \textbf{H\"{o}lder Continuity Problem}. We summarize the main results as follows.
\begin{theorem}\label{th1}
Let $d\geq 2$. Assume that $(s,p,r)$ satisfies
\begin{align}\label{cond}
  1\leq p,r\leq \infty, s>\max\{1+\frac dp,\frac32\} \quad \mathrm{or} \quad 1\leq p<2d,r=1,s=1+\frac{d}{p}.
\end{align}
For any $\alpha\in(0,1)$, there exits $u_0\in B^s_{p,r}(\R^d)$  such that the data-to-solution map $u_0\mapsto \mathbf{S}_{t}(u_0)\in \mathcal{C}([0,T];B^s_{p,r})$ of the Cauchy problem \eqref{0}
satisfies
\bbal
\limsup_{t\to0^+}\frac{\|\mathbf{S}_{t}(u_0)-u_0\|_{B^s_{p,r}}}{t^\alpha}=+\infty.
\end{align*}
\end{theorem}
For the sake of simplicity, we first transform the Euler-Poincar\'{e} equations \eqref{0} into transport type equations of the velocity $u(t,x)$. According to Yan-Yin \cite{yy}, \eqref{0} is equivalent to the nonlocal system
\begin{equation}\label{ep}
  \begin{cases}
    \partial_tu+u\cdot \nabla u= Q(u,u)+R(u,u):=\mathbf{P}(u), \qquad (t,x)\in \R^+\times \R^d,\\
  u(0,x)=u_0,\qquad x\in \R^d,
  \end{cases}
  \end{equation}
here
\bbal
&Q(u,v)=-(1-\De)^{-1}\mathrm{div}\Big(\nabla u\nabla v+\nabla u(\nabla v)^T-(\nabla u)^T\nabla v
 -(\mathrm{div} u)\nabla v+\frac12\mathbf{I}(\nabla u:\nabla v)\Big),
\\&R(u,v)=-(1-\De)^{-1}\Big((\mathrm{div} u)v+ (\na u)^T\cd v\Big).
\end{align*}
Where, we have used the following notations
\bbal
&(\na u^T)_{i,j}=\pa_{x_i}u_j, \quad (u\cd \na v)_i=\sum^d_{k=1}u_k\pa_{x_k}u_i, \quad (\na u\na v)_{ij}=\sum^d_{k=1}\pa_{x_i}u_k\pa_{x_k}v_j,
\\&\na u:\na v=\sum^d_{i,j=1}\pa_{x_i}u_j\pa_{x_i}v_j, \quad \big((\na u)^T\cd v\big)_i=\sum^d_{j=1}\pa_{x_i}u_jv_j.
\end{align*}

\section{Preliminaries}\label{sec2}

In this section, we will recall some facts about the Littlewood-Paley decomposition, the nonhomogeneous Besov spaces and their some useful properties. For more details, the readers can refer to \cite{bcd}.

There exists a couple of smooth functions $(\chi,\varphi)$ valued in $[0,1]$, such that $\chi$ is supported in the ball $\mathcal{B}\triangleq \{\xi\in\mathbb{R}^d:|\xi|\leq \frac 4 3\}$, and $\varphi$ is supported in the ring $\mathcal{C}\triangleq \{\xi\in\mathbb{R}^d:\frac 3 4\leq|\xi|\leq \frac 8 3\}$. Moreover,
$$\forall\,\ \xi\in\mathbb{R}^d,\,\ \chi(\xi)+{\sum\limits_{j\geq0}\varphi(2^{-j}\xi)}=1,$$
$$\forall\,\ 0\neq\xi\in\mathbb{R}^d,\,\ {\sum\limits_{j\in \Z}\varphi(2^{-j}\xi)}=1,$$
$$|j-j'|\geq 2\Rightarrow\textrm{Supp}\,\ \varphi(2^{-j}\cdot)\cap \textrm{Supp}\,\ \varphi(2^{-j'}\cdot)=\emptyset,$$
$$j\geq 1\Rightarrow\textrm{Supp}\,\ \chi(\cdot)\cap \textrm{Supp}\,\ \varphi(2^{-j}\cdot)=\emptyset.$$
Then, we can define the nonhomogeneous dyadic blocks $\Delta_j$ and nonhomogeneous low frequency cut-off operator $S_j$ as follows:
$$\Delta_j{u}= 0,\,\ if\,\ j\leq -2,\quad
\Delta_{-1}{u}= \chi(D)u=\mathcal{F}^{-1}(\chi \mathcal{F}u),$$
$$\Delta_j{u}= \varphi(2^{-j}D)u=\mathcal{F}^{-1}(\varphi(2^{-j}\cdot)\mathcal{F}u),\,\ if \,\ j\geq 0,$$
$$S_j{u}= {\sum\limits_{j'=-\infty}^{j-1}}\Delta_{j'}{u}.$$

\begin{definition}[\cite{bcd}]\label{de2.3}
Let $s\in\mathbb{R}$ and $1\leq p,r\leq\infty$. The nonhomogeneous Besov space $B^s_{p,r}$ consists of all tempered distribution $u$ such that
\begin{align*}
||u||_{B^s_{p,r}(\R^d)}\triangleq \Big|\Big|(2^{js}||\Delta_j{u}||_{L^p(\R^d)})_{j\in \Z}\Big|\Big|_{\ell^r(\Z)}<\infty.
\end{align*}
\end{definition}

Then, we have the following product laws.
\begin{lemma}[\cite{bcd} and Lemma 2.7, \cite{Li-Yin17}]\label{le-pro}
(1) For any $s>0$ and $1\leq p,r\leq\infty$, there exists a positive constant $C=C(d,s,p,r)$ such that
$$\|uv\|_{B^s_{p,r}(\mathbb{R}^d)}\leq C\Big(\|u\|_{L^{\infty}(\mathbb{R}^d)}\|v\|_{B^s_{p,r}(\mathbb{R}^d)}+\|v\|_{L^{\infty}(\mathbb{R}^d)}\|u\|_{B^s_{p,r}(\mathbb{R}^d)}\Big).$$
(2) Let $d\geq 2$ and  $(s,p,r)$ satisfies \eqref{cond} . Then, we have
$$||uv||_{B^{s-2}_{p,r}(\mathbb{R}^d)}\leq C||u||_{B^{s-1}_{p,r}(\mathbb{R}^d)}||v||_{B^{s-2}_{p,r}(\mathbb{R}^d)}, \qquad .$$
\end{lemma}

\begin{lemma}[see \cite{bcd}]\label{lm5}
For $1\leq p\leq \infty$ and $s>0$. There exists
a constant $C$, depending continuously on $p$ and $s$, we have
\bbal
\f\|2^{j s}\left\|[\Delta_j,v]\cd\na f\right\|_{L^{p}(\mathbb{R}^d)}\g\|_{\ell^{\infty}} \leq C\big(\|\na v\|_{L^{\infty}(\mathbb{R}^d)}\|f\|_{B_{p, \infty}^{s}(\mathbb{R}^d)}+\|\na f\|_{L^{\infty}(\mathbb{R}^d)}\|\na v\|_{B_{p, \infty}^{s-1}(\mathbb{R}^d)}\big),
\end{align*}
where we denote the standard commutator
$$[\Delta_j,v]\cd\na f=\Delta_j(v\cd \na f)-v\cd\Delta_j\na f.$$
\end{lemma}

\section{Proof of the main theorem}\label{sec3}

Let $\widehat{\phi}\in \mathcal{C}^\infty_0(\mathbb{R})$ be an even, real-valued and non-negative function on $\R$ and satisfy
\begin{numcases}{\widehat{\phi}(\xi)=}
1,&if $|\xi|\leq \frac{1}{4}$,\nonumber\\
0,&if $|\xi|\geq \frac{1}{2}$.\nonumber
\end{numcases}
Motivated by \cite{lyz1,ldl}, we
define the function $f_n(x)$ by
$$f_n(x)=\phi(x_1)\cos \f(\frac{17}{12}2^{n}x_1\g)\phi(x_2)\cdots \phi(x_d)\quad\text{with}\quad n\gg1.$$
Notice that  $\varphi(\xi)\equiv 1$ for $\frac43\leq |\xi|\leq \frac32$ and
\bbal
\mathrm{supp} \ \hat{f}_n\subset \Big\{\xi\in\R^d: \ \frac{17}{12}2^n-\frac12\leq |\xi|\leq \frac{17}{12}2^n+\frac12\Big\}\subset\Big\{\xi\in\R^d: \ \frac{4}{3}2^{n}\leq |\xi|\leq \frac{3}{2}2^{n}\Big\},
\end{align*}
we have
\begin{align}\label{cl}
{\Delta_j(f_n)=\mathcal{F}^{-1}\f(\varphi(2^{-j}\cdot)\widehat{f}_{n}\g)=}\begin{cases}
f_n, &\text{if}\; j=n,\\
0, &\text{if}\; j\neq n.
\end{cases}
\end{align}

\begin{lemma}\label{le5} Assume that $(s,p,r)$ satisfies \eqref{cond}.
Define the initial data $u_0(x)$ as
\bal\label{lyz-u0}
u_0(x):=\Big(f,0,\cdots,0\Big),
\end{align}
with
\bbal
f=\sum\limits^{\infty}_{n=3}n^{-2}2^{-ns} \phi(x_1)\cos \f(\frac{17}{12}2^{n}x_1\g)\phi(x_2)\cdots \phi(x_d).
\end{align*}
Then there exists some sufficiently large $n\in \mathbb{Z}^+$ and some sufficiently enough $c>0$  such that
\bbal
&\|u_0\|_{B^{s}_{p,r}}\leq C, \quad \|u_0\cd\na \De_{n}u_0\|_{L^p}\geq cn^{-2}2^{n(1-s)},
\end{align*}
where $C$ and $c$ are some positive constants.
\end{lemma}

{\bf Proof.}\;
According to \eqref{cl} yields
\begin{align}\label{u5}
\De_{n}f(x)&=n^{-2}2^{-ns} \phi(x_1)\cos \f(\frac{17}{12}2^{n}x_1\g)\phi(x_2)\cdots \phi(x_d).
\end{align}
By the definition of ${B}_{p,r}^{s}$ and , we deduce that
\begin{align*}
 \|u_{0}\|_{{B}_{p,r}^{s}(\R^d)}&=\f\|2^{js}\|\Delta_{j}f\|_{L^p(\R^d)}\g\|_{\ell^r(j\geq1)}
\leq \f\|\frac{1}{j^{2}}\g\|_{\ell^r(j\geq1)} \|\phi\|^d_{L^{p}(\R)}\leq C\|\phi\|^d_{L^{p}(\R)}.
\end{align*}
It is easy to show that
\begin{align*}
&u_{0}\cd \nabla \Delta_nu_{0}=\Big(f\partial_{x_1} \De_{n}f,0,\cdots,0\Big).
\end{align*}
From \eqref{u5}, we have
\begin{align*}
n^2f\pa_{x_1}\De_{n}f&=2^{-ns} f(x)\phi'(x_1)\cos \f(\frac{17}{12}2^{n}x_1\g)\phi(x_2)\cdots \phi(x_d)
\\& \quad -\frac{17}{12}2^{n}2^{-ns} f(x)\phi(x_1)\sin \f(\frac{17}{12}2^{n}x_1\g)\phi(x_2)\cdots \phi(x_d).
\end{align*}
Since $f(x)$ is a real-valued and continuous function on $\R^d$, then there exists some $\delta>0$ such that
\begin{align}\label{yh}
&|f(x)|\geq \fr{1}{2}|f(0)|=\fr{1}{2}\phi^d(0)\sum\limits^{\infty}_{n=3}n^{-2}2^{-ns}=:c_0\quad\text{ for any }  x_1\in B_{\delta}(0).\end{align}
Thus we have from \eqref{yh}
\begin{align*}
n^2\|u_0\cd\na\De_{n}u_0\|_{L^p}
&\geq c_02^{n}2^{-ns} \f\|\phi(x_1)\sin \f(\frac{17}{12}2^{n}x_1\g)\g\|_{L^p(B_{\delta}(0))}||\phi(x_2)\cdots \phi(x_d)||_{L^p(\R^{d-1})}
\\& \quad-C2^{-ns}\f\| \phi'(x_1)\cos \f(\frac{17}{12}2^{n}x_1\g)\g\|_{L^p(\R)}||\phi(x_2)\cdots \phi(x_d)||_{L^p(\R^{d-1})}\\
&\geq (c2^{n}-C)2^{-ns} .
\end{align*}
We choose $n$ large enough such that $C<\frac{c}{2}2^{n}$ and then finish the proof of Lemma \ref{le5}.

\begin{proposition}\label{pro3.1}
Assume that $u_0$ satisfies \eqref{lyz-u0}. Under the assumptions of Theorem \ref{th1}, we have
\begin{align}\label{w}
&\|\mathbf{S}_{t}\left(u_{0}\right)-u_0\|_{B^{s-1}_{p,r}}\lesssim t, \qquad \|\mathbf{w}\|_{B^{s-2}_{p,r}}\lesssim t^{2},
\end{align}
here and in what follows we denote
$$\mathbf{w}:=\mathbf{S}_{t}(u_0)-u_0-t\mathbf{\widetilde{u}}_0\quad \text{with}\quad\mathbf{\widetilde{u}}_0:=\mathbf{P}(u_0)-u_0\cd\na u_0.$$
\end{proposition}
{\bf Proof.}\; For simplicity, we denote $u(t):=\mathbf{S}_t(u_0)$ here and in what follows. Notice that $(s,p,r)$ satisfies \eqref{cond}, then using the local well-posedness result (see \cite{yy}), we know that there exists a positive time $T$ such that $u(t)\in \mathcal{C}([0,T];B_{p,r}^s)$. Furthermore, it holds that
\bbal
\|u(t)\|_{L^\infty_TB^s_{p,r}}\leq C\|u_0\|_{B^s_{p,r}}\leq C.
\end{align*}
Using the Newton-Leibniz formula and Lemma \ref{le-pro}, we obtain from \eqref{ep} that
\bal\label{s}
\|u(t)-u_0\|_{B^{s-1}_{p,r}}
&\leq \int^t_0\|\pa_\tau u\|_{B^{s-1}_{p,r}} \dd\tau
\nonumber\\&\leq \int^t_0\|\mathbf{P}(u)\|_{B^{s-1}_{p,r}} \dd\tau+ \int^t_0\|u\cd\na u\|_{B^{s-1}_{p,r}} \dd\tau
\nonumber\\&\les t\|u\|^{2}_{L_t^\infty B^{s}_{p,r}}
\les t\|u_0\|^{2}_{B^{s}_{p,r}}\les t.
\end{align}
By the Newton-Leibniz formula and Lemma \ref{le-pro} again, we have from \eqref{s} that
\begin{align}\label{p1}
\|\mathbf{w}\|_{B^{s-2}_{p,r}}
&\leq \int^t_0\|\partial_\tau u-\mathbf{v}_0\|_{B^{s-2}_{p,r}} \dd\tau \nonumber\\
&\les \int^t_0\|\mathbf{P}(u)-\mathbf{P}(u_0)\|_{B^{s-2}_{p,r}} \dd\tau +\int^t_0\|u\cd\na u-u_0\cd\na u_0\|_{B^{s-2}_{p,r}}\dd\tau
\nonumber
\\
&\les \int^t_0\|\na(u-u_0)\|_{B^{s-2}_{p,r}}\|\na(u+u_0)\|_{B^{s-1}_{p,r}} \dd\tau +\int^t_0\|u-u_0\|_{B^{s-1}_{p,r}}\|u+u_0\|_{B^{s-1}_{p,r}}\dd\tau\nonumber
\\
&\les t^2.
\end{align}
Thus, we finish the proof of Proposition \ref{pro3.1}.

{\bf Proof of Theorem \ref{th1}:} Notice that $$\mathbf{S}_{t}(u_0)-u_0=t\mathbf{\widetilde{u}}_0+\mathbf{w} \quad \text{and}\quad\mathbf{\widetilde{u}}_0=\mathbf{P}(u_0)-u_0\cd\na u_0.$$
By the triangle inequality and Propositions \ref{pro3.1}, we deduce that
\bal\label{M}
\|\mathbf{S}_{t}(u_0)-u_0\|_{B^{s}_{p,r}}
&\geq2^{{ns}}
\big\|\De_{n}\big(\mathbf{S}_{t}(u_0)-u_0\big)\big\|_{L^p}\nonumber\\
&=2^{{ns}}\big\|\De_{n}\big(t\mathbf{\widetilde{u}}_0+\mathbf{w}\big)\big\|_{L^p}
\nonumber\\&\geq t2^{{ns}}\|\De_{n}\mathbf{\widetilde{u}}_0\|_{L^p}
-2^{{2n}}2^{{n(s-2)}}
\big\|\De_{n}\mathbf{w}\big\|_{L^p}\nonumber\\
&\geq t2^{{n}s}\|\De_{n}\big(u_0\cd\na u_0\big)\|_{L^p}-
t2^{{n}s}\|\De_{n}\big(\mathbf{P}(u_0)\big)\|_{L^p}-C2^{2{n}}\|\mathbf{w}\|_{B^{s-2}_{p,\infty}}
\nonumber\\&\geq t2^{{n}s}\|u_0\cd\na\De_{n}u_0\|_{L^p}-t2^{{n}s}\|[\De_{n},u_0]\cd\na u_0\|_{L^p}\nonumber\\
&\quad-
Ct\|\mathbf{P}(u_0)\|_{B^{s}_{p,\infty}}-C2^{2{n}}t^2\nonumber\\&\geq t2^{{n}s}\|u_0\cd\na\De_{n}u_0\|_{L^p}-Ct\big\|2^{{n}s}\|[\De_{n},u_0]\cd\na u_0\|_{L^p}\big\|_{\ell^\infty}\nonumber\\
&\quad-t\|\mathbf{P}(u_0)\|_{B^{s}_{p,\infty}}-C2^{2{n}}t^2.
\end{align}
By Lemmas \ref{le-pro}-\ref{lm5}, one has $\|\mathbf{P}(u_0)\|_{B^{s}_{p,\infty}}\les \|u_0\|^2_{B^{s}_{p,r}}\les 1$
and
\bbal
\big\|2^{{n}s}\|[\De_{n},u_0]\cd\na u_0\|_{L^p}\big\|_{\ell^\infty}\les \|\na u_0\|_{L^\infty}\|u_0\|_{B^{s}_{p,\infty}}+
\|\na u_0\|_{L^\infty}\|\na u_0\|_{B^{s-1}_{p,\infty}}\les 1.
\end{align*}
Using all the above estimates and Lemma \ref{le5} to  \eqref{M}, we obtain
\bbal
\|\mathbf{S}_{t}(u_0)-u_0\|_{B^{s}_{p,r}}\geq ctn^{-2}2^{n}-Ct-C2^{2{n}}t^2,
\end{align*}
which implies
\bbal
t^{-\alpha}\|\mathbf{S}_{t}(u_0)-u_0\|_{B^{s}_{p,r}}\geq ct^{1-\alpha}n^{-2}2^{n}-Ct^{1-\alpha}-C2^{2{n}}t^{2-\alpha}.
\end{align*}
Thus, picking $t^{1-\alpha}_n=n^32^{-n}$ with large $n$, we have
\bbal
t^{-\alpha}_n\|\mathbf{S}_{t_n}(u_0)-u_0\|_{B^{s}_{p,r}}&\geq cn-Cn^32^{-n}-Cn^{6}t_n^{\alpha}
\geq \tilde{c}n.
\end{align*}
This completes the proof of Theorem \ref{th1}.

\vspace*{1em}
\section*{Acknowledgments}
M. Li was partially supported by the Jiangxi Provincial Natural Science Foundation (No.20232BAB201013 and 20212BAB201008), and China Postdoctoral Science Foundation (No. 2023M731431), and the National Natural Science Foundation of China (No.12101271).

\section*{Declarations}
\noindent\textbf{Data Availability} No data was used for the research described in the article.

\vspace*{1em}
\noindent\textbf{Conflict of interest}
The authors declare that they have no conflict of interest.

\end{document}